\documentclass{amsart}

\theoremstyle{definition}

\theoremstyle{remark}

\def\Ind{\setbox0=\hbox{$x$}\kern\wd0\hbox to 0pt{\hss$\mid$\hss} \lower.9\ht0\hbox to 0pt{\hss$\smile$\hss}\kern\wd0} 

\def\Notind{\setbox0=\hbox{$x$}\kern\wd0\hbox to 0pt{\mathchardef \nn=12854\hss$\nn$\kern1.4\wd0\hss}\hbox to 0pt{\hss$\mid$\hss}\lower.9\ht0 \hbox to 0pt{\hss$\smile$\hss}\kern\wd0}

\numberwithin{equation}{section}

\def \a {\bar a}
\def \b {\bar b}

\def \x {\bar x}

\def \d {\delta}

\def \D {\Delta}

\def \L {\Lambda}
\def \I {\mathcal I}

\def \L {\Lambda}
\def \V {\mathcal V}

\title{Corrigendum to the paper: \\
``Geometric axioms for differentially closed fields with several commuting derivations'', Journal of Algebra, Vol. 362, 2012.}
\author{Omar Le\'on S\'anchez \\
University of Waterloo}

\begin{document}

\begin{center}
{\bf CORRIGENDUM TO THE PAPER:} \\
{\bf ``GEOMETRIC AXIOMS FOR DIFFERENTIALLY CLOSED FIELDS} \\ 
{\bf WITH SEVERAL COMMUTING DERIVATIONS''} \\
Journal of Algebra, Vol. 362, pp.107-116, 2012. \\
\end{center}

\

\begin{center}
Omar Le\'on S\'anchez \\
University of Waterloo
\end{center}

\

%\maketitle

In the proof of Lemma~2.6~(2) the iteration of the map $\tau$ was not performed properly and in fact the lemma is wrong; a counterexample is given by $f=\x_1$ and $k=2$. This error does not, however, affect the geometric characterization given in Theorem~3.4 but only the attempt in Theorem~4.3 to express it as a first-order set of axioms. That attempt is incorrect; the main problem being that in general $\tau\V(f_1,\dots,f_s)\neq \V(f_1\dots,f_s,\tau f_1,\dots,\tau f_s)$. But a different, indeed simpler, set of first-order axioms, which we will now describe, does express the geometric characterization.

\

\noindent {\bf Theorem 4.3$'$.} Suppose $\D=\{\d_1,\dots,\d_m\}$ and $(K,\D\cup\{D\})$ is a characteristic zero differential field in $m+1$ commuting derivations. Then $(K,\D\cup\{D\})\models DCF_{0,m+1}$ if and only if the following hold:
\begin{enumerate}
\item $(K,\Delta)\models DCF_{0,m}$
\item Suppose $\L$ is a characteristic set of a prime $\D$-ideal of $K\{x_1,\dots,x_n\}$, $O$ is a nonempty $\D$-open subset of $\V(\L)$ disjoint from $\V(H_\L)$, and $$W\subseteq  \V(f,\tau f :\, f\in \L)$$ is a $\D$-closed set whose projection to $\V(\L)$ contains $O$. Then there exists $\a \in O$ with $(\a,D\a) \in W$.
\end{enumerate}

\

\noindent {\it Remarks.}
\begin{enumerate}
\item[(i)] Recall that $H_\L$ is the product of the separants and initials of the elements of $\L$.
\item[(ii)] Condition~(2) of 4.3$'$ is first-order expressible in the language of differential rings. Indeed, all that needs to be checked is that ``$\L=\{f_1,\dots,f_s\}$ is a characteristic set of a prime $\D$-ideal of $K\{\x\}$" is a definable property on the coefficients of $f_1,\dots,f_s$. This is done by Tressl in \S 4 of \cite{Tr} using Rosenfeld's criterion which reduces the problem to the classical problem of checking primality in polynomial rings in finitely many variables where uniform bounds are well-known. 
\item [(iii)] These axioms for $DCF_{0,m+1}$ refer to $DCF_{0,m}$. Applying the theorem to the latter we have a similar characterization of $DCF_{0,m}$ in terms of $DCF_{0,m-1}$ plus a geometric axiom, and so on, until we get to $DCF_{0,0}:=ACF_0$. That is, the theorem leads recursively to a full set of geometric axioms. Actually, it is possible to present these axioms all at once as one scheme by allowing linear combinations over the integers of the derivations (as was done in the statement of the original Theorem~4.3, for example) but we have decided for the sake of clarity to present only the relative version in this corrigendum.
\end{enumerate}

\

\begin{proof}[Proof of Theorem~4.3$'$.]
Suppose $(K,\D\cup\{D\})$ is differentially closed, and we are given  $\L$, $O\subseteq \V(\L)\setminus\V(H_\L)$, and $W\subseteq \V(f,\tau f:\, f\in\L)$ satisfying the hypotheses of~(2). By assumption $\L$ is a characteristic set of the prime $\D$-ideal
$$[\L]:H_\L ^\infty=\{f\in K\{\x\}:H_\L^\ell f\in[\L]\text{ for some }\ell\}.$$
Let $V:=\V([\L]:H_\L^\infty)$, so $V$ is an irreducible component of $\V(\L)$ and $O\subseteq V$. Let $\widehat W$ be an irreducible component of $W$ that projects $\D$-dominantly onto $V$. 

We claim that $\tau V|_O=\V(f,\tau f: f\in\L)|_O$. Recall that, by definition, $\tau V$ is $\V(f,\tau f:f\in \mathcal I(V/K))$. It is easy to see that $\V(f,\tau f: f\in\L)=\V(f,\tau f:f\in[\L])$. So, supposing that $(\bar a,\bar b)$ is a root of $f$ and $\tau f$ for all $f\in[\L]$, and $\a\in O$, we need to show that $(\a,\b)$ is a root of $\tau g$ for all $g\in\mathcal I(V/K)$. But $\mathcal I(V/K)=[\L]:H_\L^\infty$, so $H_\L^\ell g\in[\Lambda]$ for some $\ell$. We get
\begin{eqnarray*}
0&=& \tau \big(H_\L^\ell g\big)(\bar a,\bar b) \ \ \ \ \ \ \ \ \ \text{ as $H_\Lambda^\ell g\in[\L]$}\\
&=& H_\L^\ell(\bar a)\tau g(\bar a,\bar b)+g(\bar a)\tau (H_\L^\ell)(\bar a,\bar b)\\
&=& H_\L^\ell(\bar a)\tau g(\bar a,\bar b).
\end{eqnarray*}
Since $O$ is disjoint from $\V(H_\L)$ we have that $\tau g(\bar a,\bar b)=0$, as desired.

It follows that a nonempty $\D$-open subset of $\widehat W$ is contained in $\tau V$, and hence, by irreducibility, $\widehat W\subseteq \tau V$. We can now apply Theorem~3.4 (the geometric characterization of $DFC_{0,m+1}$) to $O\subseteq V$ and $\widehat W\subseteq \tau V$ to obtain $\bar a\in O$ such that $(\bar a,D\bar a)\in \widehat W\subseteq W$, as desired.

For the converse we suppose that~(2) holds and we check the geometric characterization given in Theorem~3.4. That is, given irreducible $\D$-closed sets $V\subseteq K^n$ and $W\subseteq \tau V$, with $W$ projecting $\D$-dominantly onto $V$, we need to find a point $\a\in V$ such that $(\bar a,D\bar a)\in W$.

Let $\L$ be a characteristic set of $\I(V/K)$ and let $O$ be a nonempty $\D$-open subset of $V\setminus \V(H_\L)=\V(\L)\setminus \V(H_\L)$ that is contained in the projection of $W$ (this is possible since $W$ projects $\D$-dominantly onto $V$ and $V$ is irreducible). Applying~(2) to $\L,O$ and $W$, we obtain $\bar a\in O\subseteq V$ such that $(\bar a,D\bar a)\in W$.
\end{proof}

\

The precise changes required to make the paper formally correct are:
\begin{itemize}
\item Delete 2.6~(2), 2.7, 2.8 and 3.2 (wich are false).
\item In the proof of Remark~3.3~(1) drop the reference to 3.2 and use instead the fact that if $V$ is defined over a $\D$-subfield $F\leq K$ then $\I(V/K)=\I(V/F)K\{\x\}$.
\item Replace Theorem 4.3 and its proof by the above Theorem 4.3$'$ and the proof given here.
\end{itemize}

\

%\bibliographystyle{plain}
%\bibliography{ccs}

\end{document}